\documentclass[12pt]{article}
\usepackage{amscd}
\usepackage{amsfonts}
\usepackage{amsthm}
\usepackage{amsmath}
\usepackage{latexsym}
\usepackage{pstricks}
\usepackage[all,cmtip]{xy}
\numberwithin{equation}{section}
\theoremstyle{plain}
\newtheorem{theorem}{Theorem}[section]

\newtheorem{proposition}[theorem]{Proposition}
\newtheorem{corollary}[theorem]{Corollary}
\theoremstyle{definition}

\newtheorem*{remark}{Remark}
\newcommand\Sy{\operatorname{Sym}}
\newcommand\QS{\operatorname{QSym}}
\newcommand\NS{\operatorname{NSym}}
\newcommand\Sym{\operatorname{Symm}}

\newcommand\id{\operatorname{id}}
\newcommand\im{\operatorname{im}}

\newcommand\de{\delta}
\newcommand\ep{\epsilon}

\newcommand\ka{\kappa}
\newcommand\la{\lambda}

\newcommand\si{\sigma}

\newcommand\De{\Delta}

\newcommand\ap{A_+}
\newcommand\am{A_-}

\newcommand\A{\mathcal A}
\newcommand\B{\mathcal B}
\renewcommand\P{\mathcal P}
\renewcommand\H{\mathcal H}
\newcommand\T{\mathcal T}
\newcommand\Zp{\mathbf Z^+}
\newcommand\PP{\mathfrak P}

\newcommand\<{\langle}
\renewcommand\>{\rangle}

\newlength{\BiblioSpacing}
\renewenvironment{thebibliography}[1]{%
\begin{oldthebibliography}{#1}%
   \setlength{\parskip}{\BiblioSpacing}
   \setlength{\itemsep}{\BiblioSpacing}
}%
{%
\end{oldthebibliography}%
}
\begin{document}
\title{Rooted Trees and Symmetric Functions:
Zhao's Homomorphism and the Commutative Hexagon}
\author{Michael E. Hoffman \\
\small Dept. of Mathematics\\[-0.8 ex]
\small U. S. Naval Academy, Annapolis, MD 21402 USA\\[-0.8 ex]
\small \texttt{meh@usna.edu}}
\date{\small \today \\
\small Keywords: Connes-Kreimer Hopf algebra, rooted trees, planar rooted 
trees, quasi-symmetric functions, noncommutative symmetric functions\\
\small MR Classifications:  Primary 05C05, 16W30; Secondary 81T15}
\maketitle
\begin{abstract}
Recent work in perturbative quantum field theory has led to much study
of the Connes-Kreimer Hopf algebra.  Its (graded) dual, the
Grossman-Larson Hopf algebra of rooted trees, had already been
studied by algebraists.  L. Foissy introduced a noncommutative version
of the Connes-Kreimer Hopf algebra, which turns out to be self-dual.
Using some homomorphisms defined by the author and W. Zhao, we
describe a commutative diagram that relates the aforementioned
Hopf algebras to each other and to the Hopf algebras of symmetric
functions, noncommutative symmetric functions, and quasi-symmetric
functions.
\end{abstract}
\section{Introduction}
A. Connes and D. Kreimer \cite{CK} introduced a Hopf algebra (denoted here
by $\H_K$) to study renormalization in quantum field theory.
The Hopf algebra $\H_K$ is the free commutative algebra on rooted trees, 
with a noncommutative coproduct.
Its graded dual $\H_K^*$ is isomorphic to a Hopf algebra (which we call 
$k\T$) studied earlier by R. Grossman and R. G. Larson \cite{GL}
whose elements are rooted trees with a noncommutative product and a
cocommutative coproduct.
A noncommutative version of $\H_K$, denoted here by $\H_F$, was 
introduced by L. Foissy \cite{F1,F2}:  unlike $\H_K$, it is self-dual.
In \cite{H4} the author defined a Hopf algebra $k\P$, based on planar rooted
trees in the same way $k\T$ is based on rooted trees, which is isomorphic
to $\H_F^*\cong\H_F$.  We describe these Hopf algebras in \S3 below.
\par
The author's earlier paper \cite{H4} related the Hopf algebras of the 
preceding paragraph to the well-known Hopf algebra $\Sy$ of symmetric 
functions, and its more recent extensions:  
$\QS$, the quasi-symmetric functions, and $\NS$, the noncommutative 
symmetric functions (described in \S2 below).
In \cite{H4} the author gave a pair of commutative squares that
relate the Hopf algebras named above, and demonstrated their usefulness
for certain computations.
This paper adds some important extensions to this picture.
As described in \S4, a Hopf algebra homomorphism $Z:\NS\to k\T$ due to 
W. Zhao \cite{Z2} and its dual $Z^*:\H_K\to\QS$ link the two commutative 
squares into a single diagram, which we call the commutative hexagon.
We also give an explicit characterization of $Z^*$ and deduce several
corollaries, including an easy proof of the surjectivity of $Z^*$ (hence 
the injectivity of $Z$), and a description of $Z^*$ via quasi-symmetric
generating functions of posets.
\section{Symmetric and Quasi-Symmetric Functions}
\par
Henceforth $k$ is a field of characteristic 0.
Let $\PP$ be the subalgebra of the formal power series ring 
$k[[x_1,x_2,\dots]]$
consisting of those formal power series of bounded degree, where each $x_i$
has degree 1.
An element $f\in\PP$ is called symmetric if the coefficients in $f$ of the 
monomials
\begin{equation}
\label{mon}
x_{n_1}^{i_1}x_{n_2}^{i_2}\cdots x_{n_k}^{i_k}\quad\text{and}\quad
x_1^{i_1}x_2^{i_2}\cdots x_k^{i_k}
\end{equation}
agree for any sequence of distinct positive integers $n_1,n_2,\dots,n_k$,
and quasi-symmetric if the coefficients in $f$ of the monomials (\ref{mon}) 
agree for any 
strictly increasing sequence $n_1<n_2<\dots<n_k$ of positive integers.
The sets of symmetric and quasi-symmetric formal series (by tradition
called symmetric and quasi-symmetric functions) are denoted 
$\Sy$ and $\QS$ respectively:  both are subalgebras of $\PP$, and
evidently $\Sy\subset\QS$.
\par
As a vector space, $\QS$ is generated by the monomial quasi-symmetric 
functions $M_I$, which are indexed by compositions (finite sequences)
$I=(i_1,\dots,i_k)$ of positive integers:
\[
M_I=\sum_{n_1<n_2<\dots<n_k} x_{n_1}^{i_1}x_{n_2}^{i_2}\cdots 
x_{n_k}^{i_k} 
\]
(The degree of $M_I$ is $|I|=i_1+\dots+i_k$; we set $M_{\emptyset}=1$).
Forgetting order in a composition gives a partition:  $\Sy$ has a
basis consisting of monomial symmetric functions
\[
m_{\la}=\sum_{\pi(I)=\la} M_I,
\]
where $\pi$ is the function from compositions to partitions that
forgets order.  
For example, $m_{2,1,1}=M_{(2,1,1)}+M_{(1,2,1)}+M_{(1,1,2)}$.
\par
It is well-known that $\Sy$, as an algebra, is freely generated by
several sets of symmetric functions (see, e.g., \cite{M}), in particular
(1) the elementary symmetric functions $e_k=m_{1^k}$ (where $1^k$ means 1
repeated $k$ times); (2) the complete symmetric functions 
\[
h_k=\sum_{|\la|=k} m_{\la}=\sum_{|I|=k} M_I ;
\]
and (3) the power-sum symmetric functions $p_k=m_k$.
There is a duality between the $e_k$ and the $h_k$, reflected in the
(graded) identity
\begin{equation}
\label{ideh}
(1+e_1+e_2+\cdots)(1-h_1+h_2-\cdots)=1 .
\end{equation}
\par
There is also a well-known Hopf algebra structure on $\Sy$ \cite{Gei}.
This structure can be defined by making the elementary symmetric functions
divided powers, i.e.,
\[
\De(e_k)=\sum_{i+j=k}e_i\otimes e_j .
\]
Equivalently, the $h_i$ are required to be divided powers, or the $p_i$ 
primitives.
For this Hopf algebra structure,
\begin{equation}
\label{copm}
\De(m_{\la})=\sum_{\la=\mu\cup\nu}m_{\mu}\otimes m_{\nu},
\end{equation}
where the sum is over all pairs $(\mu,\nu)$ such that $\mu\cup\nu=\la$
as multisets.  For example, $\De(m_{2,1,1})$ is 
\[
1\otimes m_{2,1,1}+m_1\otimes m_{2,1}+m_2\otimes m_{1,1}
+m_{1,1}\otimes m_2+m_{2,1}\otimes m_1+m_{2,1,1}\otimes 1 .
\]
The Hopf algebra $\Sy$ is commutative and cocommutative, so its
antipode $S$ is an algebra isomorphism with $S^2=\id$.  
In fact, as follows from (\ref{ideh}), $S(e_i)=(-1)^i h_i$.
\par
We recall a useful duality result on graded connected Hopf
algebras from \cite{H4}.
By an inner product on a graded connected Hopf algebra $\A$, we mean
a nondegenerate symmetric linear function $(\cdot,\cdot):\A\otimes\A\to k$
such that $(a,b)=0$ for homogeneous $a,b\in\A$ of different degrees.
\begin{theorem}
\label{isothm}
Let $\A, \B$ be graded connected locally finite Hopf algebras 
over $k$ which admit inner products $(\cdot,\cdot)_{\A}$ and 
$(\cdot,\cdot)_{\B}$ respectively.  
Suppose there is a degree-preserving linear map $\psi:\A\to\B$ such that, 
for all $a_1,a_2,a_3\in\A$,
\begin{enumerate}
\item[(a)]
$(a_1,a_2)_{\A}=(\psi(a_1),\psi(a_2))_{\B}$;
\item[(b)]
$(a_1a_2,a_3)_{\A}=(\psi(a_1)\otimes\psi(a_2),\De\psi(a_3))_{\B}$;
\item[(c)]
$(a_1\otimes a_2,\De(a_3))_{\A}=(\psi(a_1)\psi(a_2),\psi(a_3))_{\B}$.
\end{enumerate}
Then $\B$ is isomorphic to $\A^*$ via the pairing 
$\<b,a\>=(b,\psi(a))_{\B}$.
\end{theorem}
We note that it follows from this result that a graded connected locally
finite Hopf algebra $\A$ is self-dual provided it admits an inner product
$(\cdot,\cdot)$ such that
\[
(a_1\otimes a_2,\De(a_3))=(a_1a_2,a_3)
\]
for all $a_1,a_2,a_3\in\A$.
In particular, note that $\Sy$ has an inner product with
\[
(e_{\la},m_{\mu})=\de_{\la,\mu}
\] 
for all partitions $\la,\mu$ (where $e_{\la}$ means $e_{\la_1}e_{\la_2}\cdots$
for $\la=\la_1,\la_2,\dots$).
Then by equation (\ref{copm}),
\[
(e_{\mu}e_{\nu},m_{\la})=(e_{\mu}\otimes e_{\nu},\De(m_{\la}))=
\de_{\mu\cup\nu,\la}
\]
so $\Sy$ is self-dual.
\par
To give $\QS$ the structure of a graded connected Hopf algebra we
define the coproduct $\De$ by
\begin{equation}
\label{qcop}
\De(M_K)=\sum_{I\sqcup J=K} M_I\otimes M_J,
\end{equation}
where $I\sqcup J$ is the juxtaposition of the compositions $I$ and $J$.
While this coproduct extends that on $\Sy$, it is not cocommutative, e.g.,
\[
\De(M_{(2,1,1)})=1\otimes M_{(2,1,1)}+M_{(2)}\otimes M_{(1,1)}+
M_{(2,1)}\otimes M_{(1)}+M_{(2,1,1)}\otimes 1 .
\]
\par
Since $\QS$ is commutative but not cocommutative, it cannot be self-dual:
in fact, its dual is the Hopf algebra $\NS$ of noncommutative symmetric
functions in the sense of Gelfand \it et al. \rm\cite{Get}.  
As an algebra $\NS$ is the noncommutative polynomial algebra 
$k\< E_1, E_2,\dots\>$, with $E_i$ in degree $i$, 
and the coalgebra structure is determined by declaring the $E_i$ 
divided powers.
\begin{theorem}
\label{NSQS}
(Cf. \cite[Theorem 6.1]{Get}) The Hopf algebra $\NS$
is dual to $\QS$ via the pairing $\<E_I,M_J\>=\de_{I,J}$, where
$E_{(i_1,\dots,i_l)}=E_{i_1}E_{i_2}\cdots E_{i_l}$.
\end{theorem}
\begin{proof}
We use Theorem \ref{isothm} with $\psi:\QS\to\NS$ defined by 
$\psi(M_I)=E_I$.
Define inner products on $\NS$ and $\QS$ by $(E_I,E_J)=\de_{I,J}$
and $(M_I,M_J)=\de_{I,J}$ respectively.
Then hypothesis (a) of Theorem \ref{isothm} is immediate, 
and hypothesis (c) follows from equation (\ref{qcop}).  
For hypothesis (b), note that each nonzero contribution to
\[
(E_{i_1}\cdots E_{i_k}\otimes E_{j_1}\cdots E_{j_l},
\De(E_{n_1}\cdots E_{n_m}))
\]
comes from a splitting of $(n_1,\dots,n_m)$ with each $n_a$ the sum of a part
of $(i_1,\dots,i_k)$ (or zero) and a part of $(j_1,\dots,j_l)$ (or zero)
so that the parts stay in order; and this is exactly when a nonzero
contribution to
\[
(M_{(i_1,\dots,i_k)}M_{(j_1,\dots,j_l)},M_{(n_1,\dots,n_m)})
\]
occurs, as follows from the quasi-shuffle multiplication on $\QS$
(see \cite{H1}).
\end{proof}
There is an abelianization homomorphism $\tau:\NS\to\Sy$ sending $E_i$ 
to the elementary symmetric function $e_i$.  
Its dual $\tau^*:\Sy\to\QS$ is the inclusion $\Sy\subset\QS$.
\section{Hopf Algebras of Rooted Trees}
\par
A rooted tree is a partially ordered set $t$ (whose elements we call
vertices) such that (1) $t$ has a unique maximal vertex (the root);
and (2) for any vertex $v$, the set of vertices exceeding $v$ 
is a chain.
If a vertex $v$ covers $w$ in the partial order, we call $v$ 
the parent of $w$ and $w$ a child of $v$.
We visualize a rooted tree as a directed graph with an edge
from each vertex to each of its children:  the root (uniquely) has
no incoming edges.
Let $\T$ be the set of (finite) rooted trees, and 
$\T_n=\{t\in\T :\ |t|=n+1\}$ the set of rooted trees with $n+1$ vertices.  
There is a graded vector space
\[
k\T=\bigoplus_{n\ge 0}k\T_n 
\]
with the set of rooted trees as basis.
We denote by $\Sym(t)$ the symmetry group of the rooted tree $t$,
i.e., the group of automorphisms of $t$ as a poset (or directed graph).
%For $v\in t$ let $t_v$ be the rooted tree consisting of $v$ and its 
%descendants (with the partial order inherited from $t$).  
%If $C(v)=\{v_1,\dots,v_k\}$ is the set of children of $v$ and 
%$SG(t,v)$ the group of permutations of $C(v)$ generated by those 
%that exchange $v_i$ with $v_j$ when $t_{v_i}$ and $t_{v_j}$ are 
%isomorphic rooted trees, then
%\[
%\Sym(t)=\prod_{\text{vertices $v$ of $t$}} SG(t,v) .
%\]
\par
For any forest (i.e., monomial in rooted trees) $t_1t_2\cdots t_k$, 
there is a rooted tree $B_+(t_1t_2\cdots t_k)$ given by attaching a 
new root vertex to each of the roots of $t_1,t_2,\dots,t_k$, e.g.,
\vskip .1in
\[
B_+(\bullet\psline{*-*}(.5,-.25)(.5,.25)
\hskip .35in ) =
\psline{*-*}(.25,0)(.5,.5)
\psline{*-*}(.5,.5)(.75,0)
\psline{*-*}(.75,0)(.75,-.5)
\hskip .4in .
\]
\vskip .2in
\par\noindent
If we let $B_+(1)=\bullet\in\T_0$, then $B_+$ becomes an isomorphism
of graded vector spaces from the symmetric algebra $S(k\T)$ to $k\T$.
Here $S(k\T)$ is graded by
\[
|t_1\cdots t_k|=|t_1|+\dots+|t_k| ,
\]
where $|t|$ is the number of vertices of the rooted tree $t$. 
\par
Grossman and Larson \cite{GL} define a product $\circ$ on $k\T$ as
follows.
For rooted trees $t$ and $t'$, let $t=B_+(t_1t_2\cdots t_n)$ and $|t'|=m$.
Then $t\circ t'$ is the sum of the $m^n$ rooted trees obtained by 
attaching each of the $t_i$ to some vertex of $t'$:  if $t=\bullet$, set
$t\circ t'=t'$.  For example,
\par\noindent
\vskip .2in
\[
\psline{*-*}(.25,0)(.25,.5)
\hskip .3in \circ
\psline{*-*}(.25,0)(.5,.5)\psline{*-*}(.5,.5)(.75,0)
\hskip .4in = \hskip .1in
\psline{*-*}(.25,0)(.5,.5)\psline{*-*}(.5,.5)(.5,0)\psline{*-*}(.5,.5)(.75,0)
\hskip .4in + \quad 2
\psline{*-*}(.25,0)(.5,.5)\psline{*-*}(.5,.5)(.75,0)\psline{*-*}(.75,0)(.75,-.5)
\hskip .4in ,
\]
while
\vskip .2in
\[
\psline{*-*}(.25,0)(.5,.5)\psline{*-*}(.5,.5)(.75,0)
\hskip .4in \circ
\psline{*-*}(.25,0)(.25,.5)
\hskip .3in = \hskip .1in
\psline{*-*}(.25,0)(.5,.5)\psline{*-*}(.5,.5)(.5,0)\psline{*-*}(.5,.5)(.75,0)
\hskip .4in + \quad 2
\psline{*-*}(.25,0)(.5,.5)\psline{*-*}(.5,.5)(.75,0)\psline{*-*}(.75,0)(.75,-.5)
\hskip .5in +
\psline{*-*}(.25,-.5)(.5,0)\psline{*-*}(.5,0)(.5,.5)\psline{*-*}(.5,0)(.75,-.5)
\hskip .4in .
\]
\vskip .2in
\par\noindent
This noncommutative product makes $k\T$ a graded algebra with
two-sided unit $\bullet$.  There is a coproduct $\De$ on $k\T$ defined
by $\De(\bullet)=\bullet\otimes\bullet$ and
\begin{equation}
\label{cop}
\De(B_+(t_1t_2\cdots t_k))=
\sum_{I\cup J=\{1,2,\dots,k\}}B_+(t(I))\otimes B_+(t(J)) ,
\end{equation}
where $t_1,\dots,t_k$ are rooted trees, 
the sum is over all disjoint pairs $(I,J)$ of subsets of 
$\{1,2,\dots,k\}$ such that $I\cup J=\{1,2,\dots,k\}$, and 
$t(I)$ means the product of $t_i$ for $i\in I$ (and $t(\emptyset)=1$).
As proved in \cite{GL}, the vector space $k\T$ with product $\circ$ and
coproduct $\De$ is a graded connected Hopf algebra.
\par
It is evident from definition (\ref{cop}) that any rooted tree of the 
form $B_+(t)$, i.e., any rooted tree in which the root vertex has only one 
child, is a primitive element of the Hopf algebra $k\T$.  
Let $\P\T$ be the set of such primitive trees, graded like $\T$ by the 
number of non-root vertices.
Then the vector space of primitives $\P(k\T)$ of the Hopf algebra $k\T$
has $\P\T$ as basis, i.e., $\P(k\T)=k\P\T$:  for proof see
\cite[Theorem 4.1]{GL} or \cite[Prop. 4.2]{H2}.
\par
As a graded algebra, the Connes-Kreimer Hopf algebra $\H_K$ 
is $S(k\T)$ with the grading discussed above.
The coproduct on $\H_K$ can be described recursively by setting 
$\De(1)=1\otimes 1$ and 
\begin{equation}
\label{cocy}
\De(t)=t\otimes 1+(\id\otimes B_+)\De(B_-(t)),
\end{equation}
\par\noindent
for rooted trees $t$, where $B_-$ is the inverse of $B_+$ and it
is assumed that $\De$ acts multiplicatively on products of rooted
trees.  
Equation (\ref{cocy}) implies that the ``ladders'' $\ell_i=B_+^{i-1}(\bullet)$
are divided powers, so $\phi(e_i)=\ell_i$ defines a Hopf algebra homomorphism 
$\phi: \Sy\to k\T$.
\par
There is an inner product $(\cdot,\cdot)$ on $k\T$ given by 
\begin{equation}
\label{tip}
(t,t')=\begin{cases} |\Sym(t)|,& t=t',\\
0,&\text{otherwise.}\end{cases}
\end{equation}
Since $\Sym(B_+(t))\cong\Sym(t)$ for rooted trees $t$, $B_+$ is
an isometry for $(\cdot,\cdot)$.  In fact, we can extend
this inner product to $S(k\T)=\H_K$ by setting
\[
(u,v)=(B_+(u),B_+(v)) 
\]
for forests $u,v$.
In \cite[Theorem 4.2]{H4} (cf. \cite[Prop. 4.4]{H2})
it is shown that Theorem \ref{isothm} applies to $\H_K$ and $k\T$ 
with these inner products and $\psi=B_+$, giving the following result.
\begin{theorem}
\label{dualCKGL}
The Hopf algebra $\H_K$ is isomorphic to the dual of $k\T$ via
the pairing $\<t,u\>=(t,B_+(u))$.
\end{theorem}
A planar rooted tree is a particular realization of
a rooted tree in the plane, i.e., we consider 
\vskip .1in
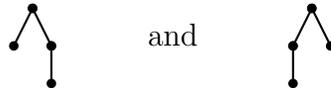
\[
\psline{*-*}(.25,0)(.5,.5)
\psline{*-*}(.5,.5)(.75,0)
\psline{*-*}(.75,0)(.75,-.5)
\hskip .8in
\text{and}
\hskip .4in
\psline{*-*}(.25,0)(.25,-.5)
\psline{*-*}(.25,0)(.5,.5)
\psline{*-*}(.5,.5)(.75,0)
\]
\vskip .2in
\par\noindent
distinct planar rooted trees.  
In parallel to $\T$, we define $\P$ to be the graded set of planar 
rooted trees, and $k\P$ the corresponding graded vector space.  
The tensor algebra $T(k\P)$ can
be regarded as the algebra of ordered forests of planar rooted trees,
and there is a linear map $B_+:T(k\P)\to k\P$ that makes a planar
rooted tree out of an ordered forest $T_1\cdots T_k$ of planar rooted 
trees by attaching a new root vertex to the root of each $T_i$.
With the same conventions about grading as above, $B_+$ is an 
isomorphism of graded vector spaces.
\par
There is an analogue of the Grossman-Larson product for planar rooted trees.  
Given planar rooted trees $T,T'$ with $T=B_+(T_1\dots T_n)$, let 
$T\circ T'$ be the sum of the 
\[
\binom{2m+n-2}{n}
\]
planar rooted trees formed by attaching $T_1,\dots,T_n$, in order,
to the vertices of $T'$.  Note that if $T'$ has $m$ vertices, then it
has a total of $2m-1$ ``attachment points'' with a natural order:  for 
example, if
\vskip .2in
\[
T'=
\hskip .2in
\psline{*-*}(-.25,0)(0,.5)
\psline{*-*}(0,.5)(.25,0)
\]
\par\noindent
then there are three attachment points $a_1,a_2,a_3$ for the root vertex,
one attachment point $b_1$ for its left child, and one attachment point
$c_1$ for its right child, with the natural order $a_1<b_1<a_2<c_1<a_3$.
Thus 
\vskip .2in
\[
\psline{*-*}(0,0)(0,.5)
\hskip .2in \circ
\hskip .2in
\psline{*-*}(-.25,0)(0,.5)
\psline{*-*}(0,.5)(.25,0)
\hskip .2in =
3\hskip .2in 
\psline{*-*}(-.25,0)(0,.5)
\psline{*-*}(0,0)(0,.5)
\psline{*-*}(.25,0)(0,.5)
\hskip .2in +
\psline{*-*}(.25,0)(.25,-.5)
\psline{*-*}(.25,0)(.5,.5)
\psline{*-*}(.5,.5)(.75,0)
\hskip .4in +
\psline{*-*}(.25,0)(.5,.5)
\psline{*-*}(.5,.5)(.75,0)
\psline{*-*}(.75,0)(.75,-.5)
\hskip .5in ,
\]
while
\vskip .2in
\[
\psline{*-*}(-.25,0)(0,.5)
\psline{*-*}(0,.5)(.25,0)
\hskip .2in \circ
\hskip .1in
\psline{*-*}(0,0)(0,.5)
\hskip .2in =
3\hskip .2in 
\psline{*-*}(-.25,0)(0,.5)
\psline{*-*}(0,0)(0,.5)
\psline{*-*}(.25,0)(0,.5)
\hskip .2in +
\psline{*-*}(.25,0)(.25,-.5)
\psline{*-*}(.25,0)(.5,.5)
\psline{*-*}(.5,.5)(.75,0)
\hskip .4in +
\psline{*-*}(.25,0)(.5,.5)
\psline{*-*}(.5,.5)(.75,0)
\psline{*-*}(.75,0)(.75,-.5)
\hskip .4in +
\hskip .2in
\psline{*-*}(-.25,-.5)(0,0)
\psline{*-*}(.25,-.5)(0,0)
\psline{*-*}(0,0)(0,.5)
\hskip .2in .
\]
\vskip .2in
\par\noindent
Now we make $k\P$ a coalgebra by defining a coproduct $\De$
on planar rooted trees by
\[
\De(T)=\sum_{i=0}^n B_+(T_1\cdots T_k)\otimes B_+(T_{k+1}\cdots T_n)
\]
where $B_-(T)=T_1\cdots T_n$.  In \cite{H4} the following result
is proved.
\begin{theorem} The product $\circ$ and coproduct $\De$ make $k\P$
a graded connected Hopf algebra.
\end{theorem}
\par
As an algebra, the Foissy Hopf algebra $\H_F$ is the tensor algebra 
$T(k\P)$.  
The coalgebra structure can be defined by the same equation 
(\ref{cocy}) as for $\H_K$,
except that rooted trees are replaced by planar rooted trees,
and the forests are ordered.  Application of Theorem \ref{isothm}
with $\psi=B_+$ and appropriate inner products gives the following
result \cite[Theorem 5.2]{H4}.
\begin{theorem} The Hopf algebra $(k\P,\circ,\De)$ is dual to $\H_F$.
\end{theorem}
On the other hand, L. Foissy constructs in \cite[\S6]{F1} an inner product 
$(\cdot,\cdot)_F$ on $\H_F$ with 
\[
(F_1F_2,F_3)_F=(F_1\otimes F_2,\De(F_3))_F
\]
for ordered forests $F_1,F_2,F_3$.  This establishes the following
result.
\begin{theorem} The Foissy Hopf algebra $\H_F$ is self-dual.
\end{theorem}
\section{The Commutative Hexagon}
The ``ladders'' $\ell_i$, considered as elements of $\H_F$, are divided
powers, so there is a Hopf algebra homomorphism $\Phi:\NS\to\H_F$ 
with $\Phi(E_i)=\ell_i$.
If $\rho:\H_F\to\H_K$ sends each planar rooted tree to the corresponding
rooted tree and forgets order in products, then $\rho\Phi=\phi\tau$.
Taking duals, we have the commutative diagram of Hopf algebras
\begin{equation}
\label{d2}
\begin{split}
\xymatrixcolsep{3pc}
\xymatrixrowsep{3pc}
\xymatrix{
\QS & k\P \ar[l]_<<<<<<{\Phi^*}\\
\Sy \ar[u]^{\tau^*} & k\T \ar[l]_<<<<<<{\phi^*} \ar[u]_{\rho^*}}
\end{split}
\end{equation}
in which the maps can be described as follows.
As noted earlier, $\tau^*$ is the inclusion $\Sy\subset\QS$.
For rooted trees $t$,
\begin{equation}
\label{Kim}
\phi^*(t)=\begin{cases} 
|\Sym(B_+(\ell_{\la}))|m_{\la},&\text{if $t=B_+(\ell_{\la})$ 
for some partition $\la$;}\\
0,&\text{otherwise;}\end{cases}
\end{equation}
where $\ell_{\la}=\ell_{\la_1}\ell_{\la_2}\cdots \ell_{\la_k}$
for a partition $\la=\la_1,\la_2,\dots,\la_k$.
On the other hand, for planar rooted trees $T$,
\[
\Phi^*(T)=\begin{cases} M_I,&\text{if $T=B_+(\ell_I)$ for some composition $I$;}\\
0,&\text{otherwise;}\end{cases}
\]
where for a composition $I=(i_1,i_2,\dots,i_k)$ we define the ordered forest
$\ell_I$ of planar rooted trees by
$\ell_I=\ell_{i_1}\ell_{i_2}\cdots \ell_{i_k}$.
For a rooted tree $t$, 
\[
\rho^*(t)=|\Sym(t)|\sum_{T\in \rho^{-1}(t)} T .
\]
\par
From (\ref{Kim}) it follows that $\phi^*$ sends the element
\[
\ka_n=\sum_{t\in\T_n}\frac{t}{|\Sym(t)|}\in k\T 
\]
to $h_n\in\Sy$.
Defining $\ep_n\in k\T$ inductively by $\ep_0=\bullet$ and
\begin{equation}
\label{eps}
\ep_n=\ka_1\circ\ep_{n-1}-\ka_2\circ\ep_{n-2}+\cdots+(-1)^{n-1}
\ka_n ,\quad n\ge 1,
\end{equation}
the identity (\ref{ideh}) implies $\phi^*(\ep_n)=e_n$.
From \cite{H4} we cite the following.
\begin{proposition}
The elements $\ka_n$ and $\ep_n$ are divided powers, and satisfy
\begin{enumerate}
\item
$\ep_n=(-1)^n S(\ka_n)$
\item
$n!\ep_n=B_+(\ell_1^n)$.
\end{enumerate}
\end{proposition}
W. Zhao defines a homomorphism $Z:\NS\to k\T$ with the
property that $Z(E_n)=\ep_n$ \cite[Theorem 4.6]{Z2}; by the 
preceding result, $Z$ sends the
noncommutative analogue $(-1)^nS(E_n)$ of the $n$th complete symmetric 
function to $\ka_n$.  In fact, $Z$ unites diagram (\ref{d2})
and its dual in a commutative hexagon.
\begin{theorem} The following diagram commutes:
\begin{equation}
\label{hex}
\begin{split}
\xymatrixcolsep{3pc}
\xymatrix{
&\H_F \ar[rd]^{\rho} \\
\NS \ar[ru]^{\Phi} \ar[rd]^{\tau} \ar[dd]_{Z} & & \H_K \ar[dd]^{Z^*} \\
&\Sy \ar[ru]^{\phi} \ar[rd]^{\tau^*} &\\
k\T \ar[ru]^{\phi^*} \ar[rd]_{\rho^*} & & \QS \\
& k\P \ar[ru]_{\Phi^*} }
\end{split}
\end{equation}
\end{theorem}
\begin{proof} Commutativity of (\ref{d2}) and its dual show the upper
and lower diamonds commute.  It remains to show that the left and right
triangles commute: since they are dual, it suffices to show $\phi^*Z=\tau$.
As the $E_i$ generate $\NS$, this follows from $\phi^*Z(E_i)=
\phi^*(\ep_i)=e_i=\tau(E_i)$.
\end{proof}
The diagram (\ref{hex}) has a symmetry about the center that exchanges each
Hopf algebra and homomorphism with its dual.  By tracing around (\ref{hex})
starting with $E_i\in\NS$ one sees, e.g., that $Z^*(\ell_i)=
M_{(1,\dots,1)}=e_i$.
Now $\rho$, $\Phi^*$ and $\phi^*$ are evidently surjective,
so $\rho^*$, $\Phi$ and $\phi$ are injective.  As we show below,
$Z^*$ is surjective, and so $Z$ is injective; cf. \cite[Theorem 5.1]{Z3}.
Thus $\NS$ is a sub-Hopf-algebra of $k\T$ (and $k\P$), while $\QS$ is
a quotient Hopf algebra of $\H_K$ (and $\H_F$).
We shall use the following characterization of $Z^*:\H_K\to\QS$.
\begin{theorem} 
\label{plus1}
For any monomial $u$ of $\H_K$, $Z^*B_+(u)=\ap Z^*(u)$,
where $\ap:\QS\to\QS$ is the linear map that sends $M_I$ to $M_{I\sqcup (1)}$.
Further, $Z^*$ is the only homomorphism of algebras from $\H_K$ to $\QS$ 
with this property.
\end{theorem}
\begin{proof} Let $\Pi:k\T\to k\P\T$ be projection onto the primitive
trees, i.e.,
\[
\Pi(t)=\begin{cases} t,&\text{if $t=B_+(t')$ for some rooted tree $t'$;}\\
0,&\text{otherwise.}\end{cases}
\]
Since $B_+^2(u)$ is a primitive tree for any monomial $u$ of $\H_K$,
we have
\[
(w,B_+^2(u))=(\Pi(w),B_+^2(u))
\] 
for any $w\in k\T$ (The inner product is that given by (\ref{tip})).
Observe that
\begin{equation}
\label{pi}
\Pi(\ep_{i_1}\circ\ep_{i_2}\circ\dots\circ\ep_{i_k})=
\begin{cases}
B_+(\ep_{i_1}\circ\ep_{i_2}\circ\dots\circ\ep_{i_{k-1}}),&\text{if $i_k=1$,}\\
0,&\text{otherwise.}\end{cases}
\end{equation}
Now let $P(E_1,E_2,\dots)\in\NS$ be any (noncommutative) polynomial
in the $E_i$.  We can write
\[
P(E_1,E_2,\dots)=\tilde P(E_1,E_2,\dots)E_1+R(E_1,E_2,\dots)
\]
where $\tilde P(E_1,E_2,\dots)$ and $R(E_1,E_2,\dots)$ are polynomials
in the $E_i$ such that every term of $R(E_1,E_2,\dots)$ ends in $E_i$
with $i>1$.
Then
\begin{align*}
\< P(E_1,E_2,\dots),Z^*B_+(u)\> &= 
\< P(\ep_1,\ep_2,\dots),B_+(u)\> \\
&= (P(\ep_1,\ep_2,\dots),B_+^2(u))\\
&= (\Pi P(\ep_1,\ep_2,\dots),B_+^2(u))\\
&= (B_+\tilde P(\ep_1,\ep_2,\dots),B_+^2(u))\\
&= (\tilde P(\ep_1,\ep_2,\dots),B_+(u))\\
&= \< \tilde P(E_1,E_2,\dots),Z^*(u)\> .
\end{align*}
Now in view of the pairing of Theorem \ref{NSQS} it is evident that
\begin{equation}
\label{apdual}
\ap^* (E_{j_1}E_{j_2}\cdots E_{j_k})=\begin{cases}
E_{j_1}E_{j_2}\cdots E_{j_{k-1}},&\text{if $j_k=1$};\\
0,&\text{otherwise,}
\end{cases}
\end{equation}
so the preceding argument shows
\[
\< P(E_1,E_2,\dots),Z^*B_+(u)\> = \<\ap^*P(E_1,E_2,\dots),Z^*(u)\>
\]
and the first part of the theorem follows.
\par
For uniqueness, note that $(\H_K,B_+)$ is the initial object in the 
category whose objects are pairs $(\A,\la)$ consisting of a commutative 
unitary algebra $\A$ and a linear map $\la:\A\to\A$, and whose morphisms
are the obvious commutative squares (see \cite[\S3]{Mo}).  
Thus there is a unique algebra homomorphism $\rho:\H_K\to\QS$ with 
$\rho B_+=\ap\rho$.
Since $Z^*$ has this property, $\rho=Z^*$.
\end{proof}
\par
We give three corollaries of Theorem \ref{plus1}.  
First we show $Z^*$ surjective.
\begin{corollary}
The homomorphism $Z^*:\H_K\to\QS$ is surjective.
\end{corollary}
\begin{proof} It suffices to show that any $M_I$ is in the image of $Z^*$.  
We proceed by induction on $|I|$, the case $|I|=1$ being immediate.  
Suppose the result holds for $M_I$ with $|I|<n$.  
From Theorem \ref{plus1} and the induction hypothesis, 
$M_{(a_1,\dots, a_{k-1},1)}\in\im Z^*$ when $a_1+\dots+a_{k-1}=n-1$.
Suppose inductively that $M_I\in\im Z^*$
when $|I|=n$ and the last part of $I$ is at most $m-1$.  Since
\begin{multline*}
M_{(1)}M_{(a_1,\dots,a_{k-1},m-1)}=M_{(1,a_1,a_2,\dots,a_{k-1},m-1)}+
M_{(a_1+1,a_2,\dots,a_{k-1},m-1)}+\dots+\\
M_{(a_1,\dots,a_{k-1},1,m-1)}+
M_{(a_1,\dots,a_{k-1},m)}+M_{(a_1,\dots,a_{k-1},m-1,1)},
\end{multline*}
it follows that $M_{(a_1,\dots,a_{k-1},m)}\in\im Z^*$ when
$a_1+\dots+a_{k-1}+m=n$.
\end{proof}
\par
Our next corollary gives a description of the homomorphism $Z^*$
via quasi-symmetric generating functions of posets.
Let $P$ be a finite poset, $\Zp$ the set of positive integers with
its usual order.
Call a function $\si:P\to\Zp$ strictly order-preserving if 
$\si(v)<\si(w)$ for any pair of elements $v<w$ of $P$, and let
$SOP(P)$ denote the set of such functions.
For any finite poset $P$, define 
\begin{equation}
\label{Kdef}
\bar K(P)=\sum_{\si\in SOP(P)} \prod_{v\in P} x_{\si(v)} .
\end{equation}
It is evident that $\bar K(P)\in\QS$ and 
\begin{equation}
\label{kmul}
\bar K(P\sqcup Q)=\bar K(P)\bar K(Q) ,
\end{equation}
where $P\sqcup Q$ means the disjoint union of posets $P$ and $Q$.
So we have a quasi-symmetric function $\bar K(t)$ for any rooted tree 
$t$, e.g.,
\[
\bar K\big(
\psline{*-*}(.25,0)(.5,.5)
\psline{*-*}(.5,.5)(.75,0)
\psline{*-*}(.75,0)(.75,-.5)
\hskip .4in \big)=\hskip .4in
\sum_{\substack{\si(a)>\si(b)\\ \si(a)>\si(c)>\si(d)}} x_{\si(a)}x_{\si(b)}
x_{\si(c)}x_{\si(d)}=
\]
%\vskip .1in
\begin{multline*}
M_{(1,1,1,1)}+M_{(1,2,1)}+M_{(1,1,1,1)}+M_{(2,1,1)}+M_{(1,1,1,1)}=\\
3M_{(1,1,1,1)}+M_{(2,1,1)}+M_{(1,2,1)} .
\end{multline*}
\begin{remark} This definition appears in \cite[\S7]{Z1}, but note
that the convention there is that the root is the minimal rather than
the maximal element of the poset corresponding to a rooted tree.
\end{remark}
Now extend $\bar K$ to $\H_K$ by defining $\bar K(t_1\cdots t_k)$ 
as $\bar K$ of the poset $t_1\sqcup\cdots\sqcup t_k$ for any 
rooted trees $t_1,\dots, t_k$.  Equation (\ref{kmul}) means that
$\bar K:\H_K\to\QS$ is an algebra homomorphism.
\begin{corollary} The algebra homomorphism $\bar K:\H_K\to\QS$
satisfies the equation $\bar K B_+(u)=\ap\bar K(u)$, and hence
must coincide with $Z^*$.
\end{corollary}
\begin{proof}  Let $u=t_1\dots t_k$ be a monomial in $\H_K$.
Thought of as a poset, $B_+(u)$ is $t_1\sqcup\dots\sqcup t_k$ with 
a maximal element adjoined.
But then it follows from definition (\ref{Kdef}) that every
term in $\bar K(B_+(u))$ consists of a factor $x_{i_1}\cdots x_{i_p}$
of $\bar K(u)$ times a factor $x_i$ corresponding to the maximal 
element, i.e. $i>i_j$ for $j=1,2,\dots,p$.  But this means that
$\bar K(B_+(u))=\ap(\bar K(u))$.
\end{proof}
\par
Finally, let $\am:\QS\to\QS$ be the adjoint of $\ap$ with respect
to the inner product on $\QS$ introduced in \S2, i.e.,
\[
\am(M_{(a_1,\dots,a_k)})=\begin{cases} M_{(a_1,\dots,a_{k-1})},
&\text{if $a_k=1$,}\\
0,&\text{otherwise.}\end{cases}
\]
Then $\am\ap=\id$.  Unlike $\ap$, $\am$ is a derivation (as can
be shown using the quasi-shuffle multiplication on $\QS$ and considering
cases).  Hence $\QS^0=\ker\am$ is a subalgebra of $\QS$ and in fact
$\QS=\QS^0[M_{(1)}]$ (see the discussion in \cite[\S2]{H3}), 
so we can think of $\am$ as $\partial/\partial M_{(1)}$.
(It's also true that $\am$ restricts to the derivation $p_1^\perp$ of $\Sy$
described in \cite[\S I.5, ex. 3]{M}.)
If $B_-$ is extended to $\H_K$ as a derivation, then we have the following.
\begin{corollary}
\label{aminus}
The homomorphism $Z^*$ satisfies $Z^*B_-=\am Z^*$.
\end{corollary}
\begin{proof} 
Since $B_-$ and $\am$ are both derivations, it suffices
to show $Z^*B_-(t)=\am Z^*(t)$ for any rooted tree $t$.  
This follows by applying $\am$ to both sides of
\[
Z^*(t)=Z^*(B_+B_-(t))=\ap Z^*B_-(t).
\]
\end{proof}
\par\noindent
In addition, the following analogue of equation (\ref{cocy}) holds:
\[
\De(u)=u\otimes 1+(\id\otimes\ap)\De(\am(u)),\quad\text{for $u\in\ap(\QS)$}.
\]
\par
Recall that the dual $\ap^*:\NS\to\NS$ of $\ap$ is given by equation 
(\ref{apdual}), and it is easy to see that $\am$ has dual $\am^*(u)=uE_1$.
The following result dualizes Theorem \ref{plus1} and Corollary \ref{aminus}.
\begin{proposition}
For $u\in\NS$,
\begin{enumerate}
\item 
$Z\ap^*(u)=B_-\Pi Z(u)$;
\item
$Z\am^*(u)=Z(u)\circ\ep_1=Z(u)\circ\ell_2$.
\end{enumerate}
\end{proposition}
\begin{proof} Part (2) is immediate:  for part (1), use equation
(\ref{pi}).
\end{proof}
This result should be compared to \cite[Prop. 4.5]{H2}:  in particular,
note the parallel between the 
operator $t\mapsto t\circ\ell_2$ of part (2) and the 
``growth operator'' $t\mapsto \ell_2\circ t$ of \cite[Prop. 4.5(1)]{H2}
(for the growth operator see also \cite[\S3]{CK}, \cite[\S2]{H2} and 
\cite[\S6.1]{H4}).

\end{document}